\documentclass[a4paper,12pt]{article}
\usepackage{amssymb}
\usepackage{color}

\newcommand{\R}{\mathbb{R}}

\newcommand{\N}{\mathbb{N}}

\newcommand{\beq}{\begin{equation} }
\newcommand{\eqq}{\end{equation} }
\newcommand{\cuad}{{\sqcap\kern-.68em\sqcup}}

\newtheorem{teo}{Theorem}[section]
\newtheorem{prop}{Proposition}[section]

\newtheorem{lemma}{Lemma}[section]

\newtheorem{remark}{Remark}[section]
\newcommand{\bremark}{\begin{remark} \em}
\newcommand{\eremark}{\end{remark} }

\def\beeq{\begin{equation}}
\def\eeq{\end{equation}}
\newcommand{\begeqaet}{\begin{eqnarray*}}
\newcommand{\eneqaet}{\end{eqnarray*}}

\newcommand{\Om}{{\Omega}}

\begin{document}
\begin{center}{\bf   Existence results of positive solutions
for nonlinear cooperative elliptic systems involving fractional Laplacian }\medskip

\bigskip

\bigskip

{ Alexander Quaas and Aliang Xia}

Departamento de  Matem\'atica,  Universidad T\'ecnica Federico Santa Mar\'{i}a

Casilla: V-110, Avda. Espa\~na 1680, Valpara\'{\i}so, Chile.

 {\sl ( alexander.quaas@usm.cl and  aliangxia@gmail.com)}
\end{center}

\bigskip

\begin{abstract}
In this article, we prove existence results of positive solutions for the following nonlinear elliptic problem  with gradient terms:
\begin{eqnarray*}
\left\{\begin{array}{l@{\quad }l}
(-\Delta)^\alpha u=f(x,u,v,\nabla u, \nabla v) &{\rm in}\,\,\Omega,\\
(-\Delta)^\alpha v=g(x,u,v,\nabla u, \nabla v) &{\rm in}\,\,\Omega,\\
u=v=0\,\,&{\rm in}\,\,\R^N\setminus\Omega,
 \end{array}
 \right.
 \end{eqnarray*}
where $(-\Delta)^\alpha$ denotes the fractional Laplacian and $ \Omega $ is a smooth bounded domain in $ \R^N $. It shown that under some assumptions on $ f $ and $ g $, the problem has at least one positive solution $(u,v)$. Our proof is based on  the classical scaling method of Gidas and Spruck and topological degree theory.
\end{abstract}
\date{}

Keywords: Fractional Laplacian; nonlinear elliptic system; gradient term; existence result.

Mathematics Subject Classification 2010: 35J60,  47G20.

\setcounter{equation}{0}
\section{ Introduction}

The problem of the existence of solutions  for  fractional elliptic systems of the form
\begin{eqnarray}\label{1.1}
\left\{\begin{array}{l@{\quad }l}
(-\Delta)^\alpha u=f(x,u,v) &{\rm in}\,\,\Omega,\\
(-\Delta)^\alpha v=g(x,u,v) &{\rm in}\,\,\Omega,\\
u=v=0\,\,&{\rm in}\,\,\R^N\setminus\Omega,
 \end{array}
 \right.
 \end{eqnarray}
where $ \Omega $ is a smooth bounded domain in $ \R^N $, $ N>2\alpha $  and $ \alpha\in(0,1) $, has been the object of intensive research during the last years.
Variational methods have been frequently used, since there is by now a remarkable collection of abstract results on the existence of critical points. But when come to the problem
under consideration is not of variational type, for example gradient terms are present, existence results of solutions are very little. In the case of scalar  equation
\begin{eqnarray*}
\left\{\begin{array}{l@{\quad }l}
(-\Delta)^\alpha u=f(x,u) &{\rm in}\,\,\Omega,\\
u=0\,\,&{\rm in}\,\,\R^N\setminus\Omega,
 \end{array}
 \right.
 \end{eqnarray*}
and
\begin{eqnarray*}
\left\{\begin{array}{l@{\quad }l}
(-\Delta)^\alpha u=f(x,u,\nabla u) &{\rm in}\,\,\Omega,\\
u=0\,\,&{\rm in}\,\,\R^N\setminus\Omega,
 \end{array}
 \right.
 \end{eqnarray*}
 have been studied in \cite{BDGQ} provided some growth condition of $f$ with respect to $u$ and $\nabla u$ is imposed. In our paper, we consider the system case and some results in \cite{FY} and \cite{YY} are extended to the fractional Laplacian case.

The fractional Laplacian $ (-\Delta)^\alpha $ is defined as
\beq\label{sin}
(-\Delta)^\alpha u(x)=C_{N,\alpha}P.V.\int_{\R^N}\frac{u(x)-u(y)}{|x-y|^{N+2\alpha}}dy\quad{\rm for}\,\,x\in \R^N.
\eqq
Here $ P.V. $ denotes the principal value of the integral, that for notational simplicity we omit in what follows and we without loss of generality take $ C_{N,\alpha}=1 $.
In additional, for some $ \sigma>0 $, suppose that $ u\in \mathcal{L}_\alpha(\R^N)\cap C^{2\alpha+\sigma} (\R^N)$ if $ 0<\alpha<1/2 $ or $ u\in \mathcal{L}_\alpha(\R^N)\cap C^{1,2\alpha+\sigma-1} (\R^N)$ if $ \alpha\ge1/2 $ where
\[
\mathcal{L}_\alpha(\R^N):=\left\{u\,\,\bigg\vert\,\,u:\R^N\rightarrow\R\,\,{\rm such\,\,that}\,\, \int_{\R^N}\frac{|u(y)|}{1+|y|^{N+2\alpha}}dy<\infty \right\},
\]
then the definition (\ref{sin}) is well defined (see Proposition 2.4 in \cite{S}).

If we separate the leading part in (\ref{1.1}), it becomes a system of the form
\begin{eqnarray}\label{1.2}
\left\{\begin{array}{l@{\quad }l}
(-\Delta)^\alpha u=a(x)u^{\beta_{11}}+b(x)v^{\beta_{12}}+h_1(x,u,v)&{\rm in}\,\,\Omega,\\
(-\Delta)^\alpha v=c(x)u^{\beta_{21}}+c(x)v^{\beta_{22}}+h_2(x,u,v)&{\rm in}\,\,\Omega,\\
u=v=0\,\,&{\rm in}\,\,\R^N\setminus\Omega,
 \end{array}
 \right.
 \end{eqnarray}
The above system  is not variational, and consequently we will use topological methods to prove the existence of positive solutions. The main difficulty
 when using topological methods lies in the need of obtaining a priori bounds. Here the a priori bounds for the solutions of (\ref{1.2}) are obtained via the so-called blow-up method (cf. Gidas and Spruck \cite{GS}).

Assuming the following conditions:

(A1) The coefficients $a,b,c,d:\bar{\Om}\rightarrow[0,\infty)$ are continous functions.

(A2) $\beta_{ij}\ge 0$ for $i,j=1,2$.

(A3) $h_1,h_2\in C(\Om, \R,\R)$, and there exists positive constants $c_1$ and $c_2$ such that
\begin{eqnarray*}
|h_1(x,s,t)|&\le& c_1(1+|s|^{\gamma_{11}}+|t|^{\gamma_{12}}),\\
|h_2(x,s,t)|&\le& c_2(1+|s|^{\gamma_{21}}+|t|^{\gamma_{22}}),
\end{eqnarray*}
where $1<\gamma_{ij}<\beta_{ij}$ for $i,j=1,2$.

We say that system (\ref{1.2}) is weakly coupled if $\beta_{11}>1$, $\beta_{22}>1$ and
\beq\label{1.3}
\beta_{12}<\frac{\beta_{22}-1}{\beta_{11}-1}\beta_{11}, \quad \beta_{21}<\frac{\beta_{11}-1}{\beta_{22}-1}\beta_{22}
\eqq
and it is strongly coupled if $\beta_{12}\beta_{21}>0$ and
\beq\label{1.4}
\beta_{11}<\frac{\beta_{21}+1}{\beta_{12}+1}\beta_{12}, \quad \beta_{22}<\frac{\beta_{12}+1}{\beta_{21}+1}\beta_{21}.
\eqq

Our  first main result is
\begin{teo}\label{t1}
Suppose (A1)-(A3) hold. 

(i) If system (\ref{1.2}) is weakly coupled and $a(x),d(x)\ge c_0>0$ for $x\in \bar{\Om}$, assume also that
\beq\label{1.5}
1<\beta_{11},\beta_{22}<\frac{N+2\alpha}{N-2\alpha},
\eqq
then system (\ref{1.2}) has at least one positive viscosity solution.

(ii) If  system (\ref{1.2}) is strongly coupled and $b(x),c(x)\ge c_0>0$ for $x\in \bar{\Om}$, assume also that
\beq\label{1.6}
\min\left\{\frac{2\alpha(\beta_{12}+1)}{\beta_{12}\beta_{21}-1},\frac{2\alpha(\beta_{21}+1)}{\beta_{12}\beta_{21}-1}\right\}>\frac{N-2\alpha}{2},
\eqq
then system (\ref{1.2}) has at least one positive viscosity solution.
\end{teo}

In this article, we will also consider problems with gradient terms. More precisely, we will study the following system
\begin{eqnarray}\label{1.7}
\left\{\begin{array}{l@{\quad }l}
(-\Delta)^\alpha u=a(x)u^{\beta_{11}}+b(x)v^{\beta_{12}}+h_1(x,u,v, \nabla u, \nabla v)&{\rm in}\,\,\Omega,\\
(-\Delta)^\alpha v=c(x)u^{\beta_{21}}+c(x)v^{\beta_{22}}+h_2(x,u,v, \nabla u, \nabla v)&{\rm in}\,\,\Omega,\\
u=v=0\,\,&{\rm in}\,\,\R^N\setminus\Omega,
 \end{array}
 \right.
 \end{eqnarray}
 We  should note that  a natural restriction in order that the gradient is meaningful for nonlocal problem is $\alpha>1/2$ (cf. \cite{BDGQ}). 
 To obtain the  a priori bounded of system (\ref{1.7}),  a serious difficulty comes when one proceeds to estimate the gradients of sequences of solutions
 that appear in the blow-up method. To handle it, we have to use some norm with  weights depending  on the distance of the boundary of domains involved, see \cite{GT, FY,BDGQ} and references therein. We define the weakly coupled and strongly coupled of system (\ref{1.7}) is the same as (\ref{1.3}) and (\ref{1.4}) respectively and also
 assume the following:
 
 (A4) $h_1,h_2\in C(\Om,\R,\R,\R^N,\R^N)$ are nonnegative, and there exists positive constants $c_1$ and $c_2$ such that
\begin{eqnarray*}
h_1(x,s,t,\xi,\eta)&\le& c_1(1+|s|^{\gamma_{11}}+|t|^{\gamma_{12}}+|\xi|^{\theta_{11}}+|\eta|^{\theta_{12}}),\\
h_2(x,s,t,\xi,\eta)&\le& c_2(1+|s|^{\gamma_{21}}+|t|^{\gamma_{22}}+|\xi|^{\theta_{21}}+|\eta|^{\theta_{22}}),
\end{eqnarray*}
where $1< \gamma_{ij}<\beta_{ij}$, $i,j=1,2$.

 (A5) If system (\ref{1.7}) is weakly coupled and $\theta_{ij}$ $(i,j=1,2)$ satisfying
 \begin{eqnarray*}
1< \theta_{11}&<&\frac{2\alpha\beta_{11}}{\beta_{11}+2\alpha-1},\quad 1<\theta_{22}<\frac{2\alpha\beta_{22}}{\beta_{22}+2\alpha-1},\\
1<\theta_{12}&<&\min\left\{2\alpha,\frac{2\alpha\beta_{11}(\beta_{22}-1)}{(\beta_{11}-1)(\beta_{22}+2\alpha-1)}\right\},\\
1< \theta_{21}&<&\min\left\{2\alpha,\frac{2\alpha\beta_{22}(\beta_{11}-1)}{(\beta_{11}+2\alpha-1)(\beta_{22}-1)}\right\}.
\end{eqnarray*}

  (A6) If system (\ref{1.7}) is strongly coupled and $\theta_{ij}$ $(i,j=1,2)$ satisfying
 \begin{eqnarray*}
1<\theta_{11}&<&\frac{2\alpha\beta_{12}(\beta_{21}+1)}{2\alpha\beta_{12}+\beta_{12}\beta_{21}+2\alpha-1},\\
 1<\theta_{22}&<&\frac{2\alpha\beta_{21}(\beta_{12}+1)}{2\alpha\beta_{21}+\beta_{12}\beta_{21}+2\alpha-1},\\
1<\theta_{12}&<&\min\left\{2\alpha,\frac{2\alpha\beta_{12}(\beta_{21}+1)}{2\alpha\beta_{21}+\beta_{12}\beta_{21}+2\alpha-1}\right\},\\
 1<\theta_{21}&<&\min\left\{2\alpha,\frac{2\alpha\beta_{21}(\beta_{12}+1)}{2\alpha\beta_{12}+\beta_{12}\beta_{21}+2\alpha-1}\right\}.
\end{eqnarray*}

\begin{remark}
We note that if $\alpha>1/2$ and $\beta_{ii}>1$, then
\[
1<\frac{2\alpha\beta_{ii}}{\beta_{ii}+2\alpha-1}<2\alpha
\]
for $i=1,2$. This implies that  $\theta_{ii} <2\alpha \,\,(i=1,2)$ in (A5). On the other hand, it permits us to assume that $\theta_{ii} >1 \,\,(i=1,2)$.
Similarly, we have $\theta_{ii} <2\alpha \,\,(i=1,2)$ in (A6) if $\alpha>1/2$ and we can also assume $\theta_{ii} >1 \,\,(i=1,2)$ in (A6)  if $\beta_{ij}\beta_{ji}>1\,\,(i,j=1,2,\,\,i\not=j)$.
\end{remark}

 We will prove that

\begin{teo}\label{t2}

Let $1/2<\alpha<1$. Suppose $(A1), (A2)$ and $(A4)$ hold. 

(i) If  system (\ref{1.7}) is weakly coupled, $a(x),d(x)\ge c_0>0$ for $x\in \bar{\Om}$ and $(A5)$ holds, then system (\ref{1.7}) has at least one positive viscosity solution if
\beq\label{1.8}
1<\beta_{11},\beta_{22}\le\frac{N}{N-2\alpha}.
\eqq

(ii) If system (\ref{1.7}) is strongly coupled, $b(x),c(x)\ge c_0>0$ for $x\in \bar{\Om}$ and $(A6)$ holds, then system (\ref{1.7}) has at least one positive viscosity solution if  
\beq\label{1.9}
\max\left\{\frac{2\alpha(\beta_{12}+1)}{\beta_{12}\beta_{21}-1},\frac{2\alpha(\beta_{21}+1)}{\beta_{12}\beta_{21}-1}\right\}\ge N-2\alpha.
\eqq

\end{teo}

We notice that we obtain the existence result of system (\ref{1.7}) in  a small range compared with (\ref{1.5}) and (\ref{1.6}) in Theorem \ref{t1}.  As mentioned before, we need to consider weighted norms which present some problems since the scaling needed near  the boundary  is not the same as in the interior. Therefore, we need to split our study into two parts : first, we obtain rough universal bounds for all solutions of (\ref{1.7}), by using the  well-known doubling lemma in \cite{PQS} and our problems are nonlocal which forces us to strengthen the subcritical hypothesis to (\ref{1.7}) and to require instead (\ref{1.8}) and (\ref{1.9}) (see Lemma \ref{l3.1} and Remark \ref{r}). After that, we reduce  the obtention of the priori bounds to an analysis near boundary.

\begin{remark}
In fact, Theorems \ref{t1} and \ref{t2} are also true if we replace the fractional Laplacian to the following more general operator
\[
(-\Delta)^\alpha_K u(x)=\int_{\R^N}\frac{2u(x)-u(x+y)-u(x-y)}{|y|^{N+2\alpha}}K(y)dy,
\]
where $\alpha\in (0,1)$ and $K$ is a measurable function defined in $\R^N$ satisfies $\lambda\le K(x)\le \Lambda$ in $\R^N$ for some constants $\lambda\le \Lambda$ and
\[
\lim_{x\rightarrow 0}K(x)=1,
\]
which is considered in \cite{BDGQ}. In particular, we notice that if set $K\equiv1$, $(-\Delta)^\alpha_K$ reduces to the fractional Laplacian.
\end{remark}

The paper is organized as follows. In Section 2 we recall some well known regularity results, convergence theorem, weighted norms and the Liouville type theorems in $\R^N$ for nonlocal systems. In Section 3, we obtain a priori bounds of systems (\ref{1.2}) and (\ref{1.7}) by the blow-up method. The existence results, Theorems \ref{t1} and \ref{t2}, are shown by topological degree theory  in Sections 4.

\setcounter{equation}{0}
\section{  Preliminaries}

The purpose of this section is to introduce some preliminaries. We start this section by recalling the following maximum principle.

\begin{prop}\label{p1}(\cite{BDGQ}, Lemma 7) Let
Let $\Omega$ be an open and bounded domain of $\R^N$, and assume that  $u\in C(\R^N)$ be a viscosity solution
of
\[(-\Delta)^\alpha u\ge0\,\, in\,\,\Omega\]
with $u\ge 0$ in $\R^N$. Then $u>0$ or $u\equiv0$ in $\Om$. 
\end{prop}

Next we give a $C^\beta$ estimate, which is a direct conclusion of Theorem 26 in \cite{CaS}.
\begin{teo}\label{2.12}
Let $\Omega$ be a regular domain. If $u\in C(\bar{\Omega})$ satisfies the inequalities
\[
\Delta^\alpha u\ge -C_0\quad and\quad\Delta^\alpha u\le C_0 \quad in \,\,\Omega,
\]
then for any $\Omega^\prime\Subset\Omega$ there exist constant $\beta>0$ such that $u\in C^\beta(\Omega^\prime)$
and
\[\|u\|_{C^\beta(\Omega^\prime)}\le C \{\sup_{\Omega}|u| + \|u\|_{L^\infty(\Omega)} + C_0\}	\]
for some constant $C > 0$ which depends on $N$.
\end{teo}

We also need the following regularity.
\begin{teo}\label{2.9} (\cite{QX}, Theorem 2.5)
Let $g$ bounded in $\R^N\setminus\Omega$ and $f\in C_{loc}^\beta(\Om)$. Suppose $u$ is a viscosity solution of
\[
(-\Delta)^\alpha u=f\quad in\,\, \Om, \quad\quad u=g\quad in \,\, \R^N\setminus\Om.
\]
Then there exists a $\gamma>0$ such that $u\in C^{2\alpha+\gamma}_{loc}(\Om)$.
\end{teo}

When $\alpha\in (1/2,1)$ the H\"older estimate for the solution can be improved to obtain  an estimate for the first derivertives. 
\begin{teo}\label{2.14} (\cite{K}, Theorem 1.2)
 Assume that $\alpha\in (1/2,1)$. Suppose $u$ is a viscosity solution of
\[
(-\Delta)^\alpha u=f\quad in\,\, \Om, 
\]
where $f\in L_{loc}^\infty(\Om)$
Then there exists a $\beta=\beta(N,\alpha)\in (0,1)$ such that $u\in C^{1,\beta}_{loc}(\Om)$. Moreover, for every ball $B_R\subset\subset \Om$ there exists a positive 
constant $C=C(N.\alpha,R)$ such that
\[
\|u\|_{C^{1,\beta}(\overline{B_{R/2}})}\le C\left(\|f\|_{L^\infty(B_R)}+\|u\|_{L^\infty(\R^N)}\right).
\]
\end{teo}

We are going to use the following convergence result for fractional Laplacian (see Lemma 5 in \cite{CaS} for integro differential equation).

\begin{teo}\label{2.13}
Let $\{u_k\}$, $k\in \N$ be a sequence of functions that are bounded in $\R^N$ and continuous in $\Omega$, $f_k$ and $f$ are continuous in  $\Omega$ such that\\
(1) $\Delta^\alpha u_k \le f_k$ in $\Omega$ in viscosity
sense.\\
(2) $u_k\rightarrow u$ locally uniformly in $\Omega$.\\
(3)$u_k\rightarrow u$ a.e. in $\R^N$.\\
(4) $f_k\rightarrow f$ locally uniformly in $\Omega$.\\
Then $\Delta^\alpha u\le f$ in $\Omega$ in viscosity
sense.
\end{teo}

Since the problems under consideration with  a right hand side which is possible singular at $\partial\Om$, we next introduce some norms which will 
help us to quantify the singularity of both right hand side and the gradient of the solutions in case $\alpha\in (1/2,1)$. We denote $d(x)=dist(x,\partial\Om)$ for $x\in \Om$. It is well know that
$d$ is Lipschitz continuous in $\Om$ with Lipschitz constant 1 and it is a $C^2$ function in a neighborhood of $\partial\Om$. We modify it outside this neighborhood to make it a $C^2$ function,
still with Lipschitz constant 1, and we extend it to be zero outside.

For $\tau\in \R$ and $u\in C(\Om)$, we define (cf. Chapter 6 in \cite{GT}) 
\[
\|u\|_0^{(\tau)}=\sup_\Om d(x)^\tau|u(x)|.
\]
When $u\in C^1(\Om)$ we also define
\beq\label{tau}
\|u\|_1^{(\tau)}=\sup_\Om \left(d(x)^\tau|u(x)|+d(x)^{\tau+1}|\nabla u(x)|\right).
\eqq
Then the following estimates are prove in \cite{BDGQ} for the Dirichlet problems.
\begin{lemma}\label{l2.1} (\cite{BDGQ}, Lemma 3)
Suppose $0<\alpha<1$. Let $f\in C(\Om)$ be such that $\|f\|_0^{(\tau)}<+\infty$ for some $\tau\in (\alpha,2\alpha)$. Then problem
\beq\label{2.4}
(-\Delta)^\alpha u=f\quad in\,\,\Om,\quad u=0\quad in\,\,\R^N\setminus\Om,
\eqq
admits a unique viscosity solution. Moreover, there exists a positive constant $C$ such that
\beq\label{2.5}
\|u\|_0^{(\tau-2\alpha)}\le C\|f\|_0^{(\tau)}.
\eqq
\end{lemma}

The next estimate concerns the gradient of the solutions of (\ref{2.4}) in case $\alpha\in (1/2,1)$.

\begin{lemma}\label{l2.2} (\cite{BDGQ}, Lemma 5)
Suppose $1/2<\alpha<1$. Let $f\in C(\Om)$ be such that $\|f\|_0^{(\tau)}<+\infty$ for some $\tau\in (\alpha,2\alpha)$. Then the unique solution
of (\ref{2.4}) verifies
\[
\|\nabla u\|_0^{(\tau-2\alpha+1)}\le C_0\left(\|f\|_0^{(\tau)}+\|u\|_0^{(\tau-2\alpha)}\right),
\]
where $C_0$ is a positive constant depends on $N$ and $\alpha$ but not on $\Om$.
\end{lemma}

The next lemma is devote to take care of the constant in (\ref{2.5}) when we consider (\ref{2.4}) in expanding domains, since in general it depends on $\Om$.
This is a crucial point for the scaling method to work properly in our setting. We take $\xi\in \partial \Om$, $\lambda>0$ and let
\[
\Om_\lambda:=\{x\in \R^N: \,\, \xi+\lambda x\in \Om\}.
\]
Observe that $d_\lambda(x):=dist(x,\partial \Om_\lambda)=\lambda^{-1}d(\xi+\lambda x)$.The following lemma show that the constant in (\ref{2.5}) for the solution of (\ref{2.4}) posed in 
$\Om_\lambda$ will depend on the domain $\Om$, but not on the dilation parameter  $\lambda$.

\begin{lemma}\label{l2.3} (\cite{BDGQ}, Lemma 6)
Suppose $0<\alpha<1$. For every $\tau\in (\alpha,2\alpha)$ and $\lambda_0>0$, there exist $C, \delta>0$
such that
\[
(-\Delta)^\alpha d_\lambda^{2\alpha-\tau}\ge Cd_\lambda^{-\tau}\quad  in \,\, (\Om_\lambda)_\delta,
\]
if $0<\lambda\le \lambda_0$. Moreover, if $u$ satisfying 
\[
(-\Delta)^\alpha u\le C_1d_\lambda^{-\tau}\quad  in \,\, \Om_\lambda,
\]
for some $C_1>0$ with $u=0$ in $\R^N\setminus \Om_\lambda$, then
\[
u(x)\le C_2(C_1+\|u\|_{L^\infty(\Om_\lambda)})d_\lambda^{2\alpha-\tau} \quad for\,\,x\in (\Om_\lambda)_\delta,
\]
for some $C_2>0$ depending on $\alpha,\delta,\tau $ and $C_0$.
\end{lemma}

We finish this section by listing the well-known 
Liouville type theorems of the limit systems of (\ref{1.1}) in the whole space has been considered in our previous  article \cite{QX1}. In \cite{QX1}, we proved that
\begin{teo}\label{t4}(\cite{QX1}, Theorem 1.1)
Let $p,q>0$ and $pq>1$. Suppose
\[
\frac{2\alpha(p+1)}{pq-1},\frac{2\alpha(q+1)}{pq-1}\in\left[ \frac{N-2\alpha}{2},\,\, N-2\alpha\right)
\]
and
\[
 \left(\frac{2\alpha(p+1)}{pq-1},\frac{2\alpha(q+1)}{pq-1}\right)\not=\left( \frac{N-2\alpha}{2},\,\, \frac{N-2\alpha}{2}\right).
\]
Then, for some $ \sigma>0 $, there exists no positive $ \mathcal{L}_\alpha(\R^N)\cap C^{2\alpha+\sigma} (\R^N)$ if $ 0<\alpha<1/2 $ or  in $\mathcal{L}_\alpha(\R^N)\cap C^{1,2\alpha+\sigma-1} (\R^N)$ if $ \alpha\ge1/2 $ type solution to system 
\begin{eqnarray}\label{ff}
\left\{\begin{array}{l@{\quad }l}
(-\Delta)^\alpha u=v^q &{ in}\,\,\R^N,\\
(-\Delta)^\alpha v=u^p &{ in}\,\,\R^N.
 \end{array}
 \right.
 \end{eqnarray}
\end{teo}

In this theorem, we consider the nonexistence of solutions just  in  the subregion
\[
\max\left\{\frac{2\alpha(p+1)}{pq-1},\frac{2\alpha(q+1)}{pq-1}\right\}<N-2\alpha.
\]
In order to prove Theorem \ref{t1} completely, we also need to study the nonexistence results in the following region
\beq\label{1.11}
\max\left\{\frac{2\alpha(p+1)}{pq-1},\frac{2\alpha(q+1)}{pq-1}\right\}\ge N-2\alpha,
\eqq
Using the fundamental solutions of fractional Laplacian and the comparison principle, 
we can show there are no positive viscosity supersolutions to system (\ref{1.6}) if and only if $(p,q)$ verifies (\ref{1.11}). In fact, we prove that

\begin{teo}\label{t5}
Suppose $p,q>0$ and $ pq>1$.  Then there are no positive viscosity supersolutions to system (\ref{ff}) if and only if $(p,q)$ verifies (\ref{1.11}).
\end{teo}

{\bf Proof.} 
Here we omit the proof since in very similar as the proof  of Theorem 1.3 (I) in \cite{QS1} by using the fundamental solutions of fractional Laplacian (see
Theorem 1.3 in \cite{FQ2} and also Theorem 1.3 in \cite{FQ1})  and the comparison principle (cf. Theorem 2.1 in \cite{FQ2}). $\Box$

\setcounter{equation}{0}
\section{ A priori bounds}

In this section is devote to get the a priori bounds of systems (\ref{1.2}) and (\ref{1.7}) by a blow-up method.

\begin{teo}\label{t4.1} Under the hypotheses of Theorem \ref{t1}, then each couple of positive viscosity solution $(u,v)$ of (\ref{1.2}) is bounded in the $L^\infty$-norm
by a constant.
\end{teo}

{\bf Proof.} Assume on the contrary that there exists a sequence $(u_n,v_n)$ of positive solutions of (\ref{1.2}) such that 
\[
\max\{\|u_n\|_{L^\infty(\Om)},\|v_n\|_{L^\infty(\Om)}\}\rightarrow+\infty
\]
as $n\rightarrow+\infty$. We may assume
\[
\lambda_n=\|u_n\|_{L^\infty(\Om)}^{-1/\sigma_1},
\]
if $\|u_n\|_{L^\infty(\Om)}^{\sigma_2}\ge \|v_n\|_{L^\infty(\Om)}^{\sigma_1}$ (up to a sequence), and
\[
\lambda_n=\|v_n\|_{L^\infty(\Om)}^{-1/\sigma_2}
\]
otherwise, for some constants $\sigma_1, \sigma_2>0$ which are to be determined later. Without loss of generality, we suppose that we are in the first of these two situations. 

Note that we have $\lambda_n\rightarrow 0$ as $n\rightarrow \infty$. Let $x_n\in \Om$ be a point where $u_n$ assumes its maximum.
The functions
\[
\tilde{u}_n(x)=\lambda_n^{\sigma_1} u_n(\lambda_n x+x_n) \quad{\rm and}\quad \tilde{v}_n(x)=\lambda_n^{\sigma_2} v_n(\lambda_n x+x_n)
\]
are such that $\tilde{u}_n(0)=1$ and $0\le \tilde{u}_n, \tilde{v}_n\le 1$ in $\Om$.
One also verifies that the functions $\tilde{u}_n$ and $\tilde{v}_n$ satisying
\begin{eqnarray}\label{ee}
\left\{\begin{array}{l@{\quad }l}
(-\Delta)^\alpha \tilde{u}_n=\lambda_n^{\sigma_1+2\alpha-\sigma_1\beta_{11}}a(z)\tilde{u}_n^{\beta_{11}}+\lambda_n^{\sigma_1+2\alpha-\sigma_2\beta_{12}}b(z)\tilde{v}_n^{\beta_{12}}\\
\quad\quad\quad\quad\quad\quad+\lambda_n^{\sigma_1+2\alpha}h_1(z, \lambda_n^{-\sigma_1}\tilde{u}_n,\lambda_n^{-\sigma_2}\tilde{v}_n)&{\rm in}\,\,\Omega_n,\\
(-\Delta)^\alpha \tilde{v}_n=\lambda_n^{\sigma_2+2\alpha-\sigma_1\beta_{21}}c(z)\tilde{u}_n^{\beta_{21}}+\lambda_n^{\sigma_2+2\alpha-\sigma_2\beta_{22}}b(z)\tilde{v}_n^{\beta_{22}}\\
\quad\quad\quad\quad\quad\quad+\lambda_n^{\sigma_2+2\alpha}h_2(z, \lambda_n^{-\sigma_1}\tilde{u}_n,\lambda_n^{-\sigma_2}\tilde{v}_n) &{\rm in}\,\,\Omega_n,
 \end{array}
 \right.
 \end{eqnarray}
where $z=\lambda_n x+x_n$, and
\[
\Om_n=\{x\in \R^N :\,\, x_n+\lambda_n x\in \Om\}.
\]
By assumption (A3) we have
\begin{eqnarray*}
|\lambda_n^{\sigma_1+2\alpha}h_1(z, \lambda_n^{-\sigma_1}\tilde{u}_n,\lambda_n^{-\sigma_2}\tilde{v}_n)|
\le C(\lambda_n^{\sigma_1+2\alpha}+\lambda_n^{\sigma_1+2\alpha-\sigma_1\gamma_{11}}|\tilde{u}_n|^{\gamma_{11}}+\lambda_n^{\sigma_1+2\alpha-\sigma_2\gamma_{12}}|\tilde{v}_n|^{\gamma_{12}})
\end{eqnarray*}
and a similar estimate for $\lambda_n^{\sigma_2+2\alpha}h_2$.

{\bf Case I: weakly coupled.} Choosing 
\[
\sigma_1=\frac{2\alpha}{\beta_{11}-1}\quad{\rm and}\quad \sigma_2=\frac{2\alpha}{\beta_{22}-1},
\]
since (\ref{1.2}) is a weakly coupled system, we obtain
\begin{eqnarray}\nonumber
\sigma_1+2\alpha-\sigma_1\beta_{11}&=&0, \quad\sigma_1+2\alpha-\sigma_2\beta_{12}>0,\\\label{b}
\sigma_2+2\alpha-\sigma_1\beta_{21}&>&0, \quad\sigma_2+2\alpha-\sigma_2\beta_{22}=0.
\end{eqnarray}
By assumption $\gamma_{ij}<\beta_{ij}$ and (\ref{b}), we also obtain inequalities
\begin{eqnarray*}
\sigma_1+2\alpha-\sigma_1\gamma_{11}&>&0, \quad\sigma_1+2\alpha-\sigma_2\gamma_{12}>0,\\
\sigma_2+2\alpha-\sigma_1\gamma_{21}&>&0, \quad\sigma_2+2\alpha-\sigma_2\gamma_{22}>0.
\end{eqnarray*}

For $x\in \Om$, we denote $d(x)={\rm dist}(x,\partial\Om)$. By passing to subsequences, two situations may arise: either $d(x_n)/\lambda_n\rightarrow +\infty$ or
 $d(x_n)/\lambda_n\rightarrow d\ge0$.

We first suppose the case $d(x_n)/\lambda_n\rightarrow +\infty$ holds. Then $\Om_n\rightarrow \R^N$ as $n\rightarrow +\infty$.  Since $\tilde{u}_n$ and $\tilde{v}_n$ are uniformly bounded, by Theorem \ref{2.12} with an application of Ascoli-Arzel\'a theorem and a diagonal argument, we have
$\tilde{u}_n\rightarrow u$ and $\tilde{v}_n\rightarrow v$ locally uniformly in $\R^N$. Then passing to the limit (use Theorem \ref{2.13}), we see that $(u,v)$ solves
\begin{eqnarray*}
\left\{\begin{array}{l@{\quad }l}
(-\Delta)^\alpha u=a(x_0)u^{\beta_{11}} &{\rm in}\,\,\R^N,\\
(-\Delta)^\alpha v=d(x_0)v^{\beta_{22}} &{\rm in}\,\,\R^N
 \end{array}
 \right.
 \end{eqnarray*}
in the viscosity sense. However, by Theorem 1.2 in \cite{QX} (see also \cite{CLO}),  we know this problem has no positive viscosity solutions if 
\[
1< \beta_{11},\beta_{22}<\frac{N+2\alpha}{ N-2\alpha}.
\]

If the case  $d(x_n)/\lambda_n\rightarrow d\ge0$ holds, then we may assume $x_n\rightarrow x_0\in \partial \Om$.
Without loss of generality, we may assume the $\nu(x_0)=-e_N$. In this case we consider functions
\[
\bar{u}_n(x)=\lambda_n^{\sigma_1} u_n(\lambda_n x+\xi_n) \quad{\rm and}\quad \bar{v}_n(x)=\lambda_n^{\sigma_2} v_n(\lambda_n x+\xi_n)\quad {\rm in}\,\, D_n,
\]
where $\xi_n\in\partial\Om$ is  the projection of $x_n$ on $\partial\Om $ and 
\[
D_n=\{x\in \R^N: \xi_n+\lambda_n x\in \Om \}.
\]
Observe that 
\beq\label{eee}
0\in \partial D_n,
\eqq
and 
\[
D_n\rightarrow \R^N_+=\{x\in \R^n:\,\, x_N>0\}\quad{\rm as}\,\, n\rightarrow+\infty.
\]
It also follows that $(\bar{u}_n,\bar{v}_n)$ satisfying (\ref{ee}) in $D_n$ with a  slightly different functions $h_1$ and $h_2$, but with same bounds.

Furthermore, let
\[
\tilde{x}_n=\frac{x_n-\xi_n}{\lambda_n},
\]
 we have $|\tilde{x}_n|=d(x_n)/\lambda_n$ and $\bar{u}_n(\tilde{x}_n)=1$. We claim that
\[
d:=\lim_{n\rightarrow +\infty}\frac{d(x_n)}{\lambda_n}>0.
\]
This is in particular guarantees that by passing to a subsequence $x_n\rightarrow x_0$, where $|x_0|=d>0$. Therefore, $x_0$ is in the interior of the half space $\R^N_+$.
Now we are in the position to prove the claim. Observe that by (\ref{ee}), we have
\[
(-\Delta)^\alpha \bar{u}_n\le C\le \tilde{C}d_n^{-\theta}\quad{\rm in}\,\, D_n,
\]
where $\theta\in (\alpha,2\alpha)$ and $d_n={\rm dist}(x,\partial D_n)$. By Lemma \ref{l2.3}, for fixed $\theta$, there exists a constant $C_0>0$ and $\delta>0$
such that $\bar{u}_n(x)\le C_0d_n(x)^{2\alpha-\theta}$ if $d_n(x)<\delta$. In particular, by (\ref{eee}), we know $|\tilde{x}_n|\ge d_n(\tilde{x}_n)$. Therefore, if $ d_n(\tilde{x}_n)<\delta$,
then $1=\bar{u}_n(\tilde{x}_n)\le C_0d_n(\tilde{x}_n)^{2\alpha-\theta}\le C_0|\tilde{x}_n|^{2\alpha-\theta}$. This implies $|\tilde{x}_n|$ is bounded from blow and thus $d>0$.

Now we can employ regularity Theorem \ref{2.13} as before to obtain that
$\bar{u}_n\rightarrow u$ and $\bar{v}_n\rightarrow v$ on compact sets of  $\R^N_+$, where $(u,v)$ verifies that $0\le u,v\le 1$ in $\R^N_+$, $u(x_0)=1$ and solves
\begin{eqnarray*}
\left\{\begin{array}{l@{\quad }l}
(-\Delta)^\alpha u=a(x_0)u^{\beta_{11}} &{\rm in}\,\,\R^N_+,\\
(-\Delta)^\alpha v=d(x_0)v^{\beta_{22}} &{\rm in}\,\,\R^N_+,\\
u=v=0&{\rm in}\,\,\R^N\setminus\R^N_+,
 \end{array}
 \right.
 \end{eqnarray*}
in viscosity sense. However, by Theorems 1.1 in \cite{QX}, we know the above problem have no positive viscosity solutions
if 
\[
1<\beta_{11},\beta_{22}<\frac{(N-1)+2\alpha}{ (N-1)-2\alpha}.
\]
This contradicts with our assumption since $\frac{(N-1)+2\alpha}{ (N-1)-2\alpha}>\frac{N+2\alpha}{ N-2\alpha}$ if $N>1+2\alpha$. 

{\bf Case II: strongly coupled.} Choosing 
\[
\sigma_{1}=\frac{2\alpha(\beta_{12}+1)}{\beta_{12}\beta_{21}-1}\quad{\rm and}\quad \sigma_{2}=\frac{2\alpha(\beta_{21}+1)}{\beta_{12}\beta_{21}-1}.
\]
Since the system (\ref{1.2}) is strongly coupled, we obtain
\begin{eqnarray}\nonumber
\sigma_1+2\alpha-\sigma_1\beta_{11}&>&0, \quad\sigma_1+2\alpha-\sigma_2\beta_{12}=0,\\\label{c}
\sigma_2+2\alpha-\sigma_1\beta_{21}&=&0, \quad\sigma_2+2\alpha-\sigma_2\beta_{22}>0.
\end{eqnarray}
and assumption $\gamma_{ij}<\beta_{ij}$  implies
\begin{eqnarray*}
\sigma_1+2\alpha-\sigma_1\gamma_{11}&>&0, \quad\sigma_1+2\alpha-\sigma_2\gamma_{12}>0,\\
\sigma_2+2\alpha-\sigma_1\gamma_{21}&>&0, \quad\sigma_2+2\alpha-\sigma_2\gamma_{22}>0.
\end{eqnarray*}
By a similar argument as before, we know the limit system of (\ref{ee}) is
\begin{eqnarray*}
\left\{\begin{array}{l@{\quad }l}
(-\Delta)^\alpha u=b(x_0)v^{\beta_{12}}\\
(-\Delta)^\alpha v=c(x_0)v^{\beta_{21}} 
 \end{array}
 \right.
 \end{eqnarray*}
in $\R^N$ or $\R^N_+$ with $u=v=0$ in $\R^N\setminus\R^N_+$.
Therefore, by Liouville type results in $\R^N$ (see Theorems \ref{t4} and \ref{t5}) and $\R^N_+$ (see Theorem 4.2 in \cite{QX}) and regularity results (cf. \cite{S}), we come to a contradiction
as before if (\ref{1.6}) holds. We complete the prove. $\Box$\\

Next, we prove a priori bound for system (\ref{1.7}).  Firstly, we obtain rough bounds for all solutions of the system (\ref{1.7}) which are universal, in the spirit of \cite{PQS}.

\begin{lemma}\label{l3.1}
Under assumptions in Theorem \ref{t2},
assume that positive function $u,v\in C^1(\Om)\cap L^\infty(\R^N)$ satisfying
\begin{eqnarray}\label{fe}
\left\{\begin{array}{l@{\quad }l}
(-\Delta)^\alpha u=a(x)u^{\beta_{11}}+b(x)v^{\beta_{12}}+h_1(x,u,v, \nabla u, \nabla v)&{\rm in}\,\,\Omega,\\
(-\Delta)^\alpha v=c(x)u^{\beta_{21}}+c(x)v^{\beta_{22}}+h_2(x,u,v, \nabla u, \nabla v)&{\rm in}\,\,\Omega,
 \end{array}
 \right.
 \end{eqnarray}
 in the viscosity sense, then there exists a positive constant $C$ such that
 \[
 u(x)\le C(1+{\rm dist}(x,\partial\Om)^{-\sigma_1}),\quad |\nabla u(x)|\le C(1+{\rm dist}(x,\partial\Om)^{-(\sigma_1+1)}),
 \]
 \[
 v(x)\le C(1+{\rm dist}(x,\partial\Om)^{-\sigma_2}),\quad |\nabla v(x)|\le C(1+{\rm dist}(x,\partial\Om)^{-(\sigma_2+1)}),
 \]
 for $x\in \Om$, where 
 \[
\sigma_1=\frac{2\alpha}{\beta_{11}-1}\quad{\rm and}\quad \sigma_2=\frac{2\alpha}{\beta_{22}-1},
\]
verifies (\ref{1.8}) in weakly coupled case (see (\ref{1.3})), and
\[
\sigma_{1}=\frac{2\alpha(\beta_{12}+1)}{\beta_{12}\beta_{21}-1}\quad{\rm and}\quad \sigma_{2}=\frac{2\alpha(\beta_{21}+1)}{\beta_{12}\beta_{21}-1}.
\]
verifies (\ref{1.9}) in strongly coupled case (see (\ref{1.4})).
\end{lemma}

{\bf Proof.} Assume that the Lemma fails. Then, there exist sequences of positive function $u_n,v_n\in C^1(\Om)\cap L^\infty(\R^N)$ and $y_n\in \Om$
satisfying 
\begin{eqnarray*}
\left\{\begin{array}{l@{\quad }l}
(-\Delta)^\alpha u_n=a(x)u_n^{\beta_{11}}+b(x)v_n^{\beta_{12}}+h_1(x,u_n,v_n, \nabla u_n, \nabla v_n) &{\rm in}\,\,\Om,\\
(-\Delta)^\alpha v_n=c(x)u_n^{\beta_{21}}+c(x)v_n^{\beta_{22}}+h_2(x,u_n,v_n, \nabla u_n, \nabla v_n)&{\rm in}\,\,\Om
 \end{array}
 \right.
 \end{eqnarray*}
and 
\[
M_n:=u_n^{\frac{1}{\beta_1}}+|\nabla u_n|^{\frac{1}{\beta_1+1}}+v_n^{\frac{1}{\beta_2}}+|\nabla v_n|^{\frac{1}{\beta_2+1}}
\]
satisfies
\beq\label{3.1}
M_n(y_n)>2n(1+{\rm dist}(y_n,\partial\Om)^{-1}).
\eqq
By Lemma 5.1 in \cite{PQS}, there exists a sequence of points $x_n\in \Om$ with the property that $M_n(x_n)\ge M_n(y_n)$, $M_n(x_n)>2n{\rm dist}(x_n,\partial\Om)^{-1}$
and 
\beq\label{3.2}
M_n(z)\le 2M_n(x_n) \quad {\rm in }\,\, B(x_n,nM_n(x_n)^{-1}).
\eqq
Observe that (\ref{3.1}) implies $M_n(x_n)\rightarrow +\infty$. Let $\lambda_n=M_n(x_n)^{-1}\rightarrow 0$ and define
\[
\tilde{u}_n(x)=\lambda_n^{\sigma_1}u_n(x_n+\lambda_nx)\quad{\rm and}\quad \tilde{v}_n(x)=\lambda_n^{\sigma_2}v_n(x_n+\lambda_nx)
\]
in $B_n:=\{x\in \R^N:\,\, |x|<n\}$.
Then function $(\tilde{u}_n, \tilde{v}_n)$ satisfies
\begin{eqnarray}\label{ue}
\left\{\begin{array}{l@{\quad }l}
(-\Delta)^\alpha \tilde{u}_n=\lambda_n^{\sigma_1+2\alpha-\sigma_1\beta_{11}}a(z)\tilde{u}_n^{\beta_{11}}+\lambda_n^{\sigma_1+2\alpha-\sigma_2\beta_{12}}b(z)\tilde{v}_n^{\beta_{12}}\\
\quad\quad+\lambda_n^{\sigma_1+2\alpha}h_1(z, \lambda_n^{-\sigma_1}\tilde{u}_n,\lambda_n^{-\sigma_2}\tilde{v}_n,\lambda_n^{-(\sigma_1+1)}\nabla\tilde{u}_n,\lambda_n^{-(\sigma_2+1)}\nabla\tilde{v}_n)&{\rm in}\,\,B_n,\\
(-\Delta)^\alpha \tilde{v}_n=\lambda_n^{\sigma_2+2\alpha-\sigma_1\beta_{21}}c(z)\tilde{u}_n^{\beta_{21}}+\lambda_n^{\sigma_2+2\alpha-\sigma_2\beta_{22}}b(z)\tilde{v}_n^{\beta_{22}}\\
\quad\quad+\lambda_n^{\sigma_2+2\alpha}h_2(z, \lambda_n^{-\sigma_1}\tilde{u}_n,\lambda_n^{-\sigma_2}\tilde{v}_n,\lambda_n^{-(\sigma_1+1)}\nabla\tilde{u}_n,\lambda_n^{-(\sigma_2+1)}\nabla\tilde{v}_n) &{\rm in}\,\,B_n,
 \end{array}
 \right.
 \end{eqnarray}
where $z=\lambda_n x+x_n$.
By assumption (A4) we have
\begin{eqnarray*}
|\lambda_n^{\sigma_1+2\alpha}h_1(z, \lambda_n^{-\sigma_1}\tilde{u}_n,\lambda_n^{-\sigma_2}\tilde{v}_n,\lambda_n^{-(\sigma_1+1)}\nabla\tilde{u}_n,\lambda_n^{-(\sigma_2+1)}\nabla\tilde{v}_n)|\\
\le C(\lambda_n^{\sigma_1+2\alpha}+\lambda_n^{\sigma_1+2\alpha-\sigma_1\gamma_{11}}|\tilde{u}_n|^{\gamma_{11}}+\lambda_n^{\sigma_1+2\alpha-\sigma_2\gamma_{12}}|\tilde{v}_n|^{\gamma_{12}}\\
+\lambda_n^{\sigma_1+2\alpha-(\sigma_1+1)\theta_{11}}|\nabla\tilde{u}_n|^{\theta_{11}}+\lambda_n^{\sigma_1+2\alpha-(\sigma_2+1)\theta_{12}}|\nabla\tilde{v}_n|^{\theta_{12}})
\end{eqnarray*}
and a similar estimate for $\lambda_n^{\sigma_2+2\alpha}h_2$.

{\bf Case I: weakly coupled.} Choosing 
\[
\sigma_1=\frac{2\alpha}{\beta_{11}-1}\quad{\rm and}\quad \sigma_2=\frac{2\alpha}{\beta_{22}-1},
\]
since (\ref{fe}) is a weakly coupled system, we obtain
\begin{eqnarray*}
\sigma_1+2\alpha-\sigma_1\beta_{11}&=&0, \quad\sigma_1+2\alpha-\sigma_2\beta_{12}>0,\\
\sigma_2+2\alpha-\sigma_1\beta_{21}&>&0, \quad\sigma_2+2\alpha-\sigma_2\beta_{22}=0.
\end{eqnarray*}
By assumption $\gamma_{ij}<\beta_{ij}$, we also obtain inequalities
\begin{eqnarray*}
\sigma_1+2\alpha-\sigma_1\gamma_{11}&>&0, \quad\sigma_1+2\alpha-\sigma_2\gamma_{12}>0,\\
\sigma_2+2\alpha-\sigma_1\gamma_{21}&>&0, \quad\sigma_2+2\alpha-\sigma_2\gamma_{22}>0.
\end{eqnarray*}
Using (A5), we have
\begin{eqnarray*}
\sigma_1+2\alpha-(\sigma_1+1)\theta_{11}&>&0, \quad\sigma_1+2\alpha-(\sigma_2+1)\theta_{12}>0,\\
\sigma_2+2\alpha-(\sigma_1+1)\theta_{21}&>&0, \quad\sigma_2+2\alpha-(\sigma_2+1)\theta_{22}>0.
\end{eqnarray*}
Moreover, by (\ref{3.2}),
\[
u_n(x)^{\frac{1}{\beta_1}}+|\nabla u_n(x)|^{\frac{1}{\beta_1+1}}+v_n(x)^{\frac{1}{\beta_2}}+|\nabla v_n(x)|^{\frac{1}{\beta_2+1}}\le 2, \quad x\in B_n.
\]
We also know that
\[
u_n(0)^{\frac{1}{\beta_1}}+|\nabla u_n(0)|^{\frac{1}{\beta_1+1}}+v_n(0)^{\frac{1}{\beta_2}}+|\nabla v_n(0)|^{\frac{1}{\beta_2+1}}=1.
\]
Since $\lambda_n\rightarrow 0$, $\tilde{u}_n$, $\tilde{v}_n$, $|\nabla\tilde{u}_n|$ and $|\nabla\tilde{v}_n|$ are uniformly bounded in $B_n$, by Theorem \ref{2.14}  to obtain, with the help of Ascoli-Arzel\'a's theorem and a diagonal argument,
that there exists a subsequence, still denoted $(\tilde{u}_n,\tilde{v}_n)$ such that $\tilde{u}_n\rightarrow u$ and $\tilde{v}_n\rightarrow v$ in $C^1_{loc}(\R^N)$ as $n\rightarrow +\infty$.
Since $u(0)^{\frac{1}{\beta_1}}+|\nabla u(0)|^{\frac{1}{\beta_1+1}}+v(0)^{\frac{1}{\beta_2}}+|\nabla v(0)|^{\frac{1}{\beta_2+1}}=1$, then $(u,v)\not=(0,0)$.

Next, let $(\bar{u}_n, \bar{v}_n)$ be the functions obtained by extending $(\tilde{u}_n, \tilde{v}_n)$ to be zero outsider $B_n$. Then we can check that $(\bar{u}_n, \bar{v}_n)$ satisfying
\begin{eqnarray*}
\left\{\begin{array}{l@{\quad }l}
(-\Delta)^\alpha \bar{u}_n\ge a(z)\bar{u}_n^{\beta_{11}} &{\rm in}\,\,B_n,\\
(-\Delta)^\alpha \bar{v}_n\ge d(z)\bar{v}_n^{\beta_{22}} &{\rm in}\,\,B_n.
 \end{array}
 \right.
 \end{eqnarray*}
Passing the limit by using Theorem \ref{2.13}, we have
\begin{eqnarray*}
\left\{\begin{array}{l@{\quad }l}
(-\Delta)^\alpha u\ge a(x_0)u^{\beta_{11}} &{\rm in}\,\,\R^N,\\
(-\Delta)^\alpha v\ge d(x_0)v^{\beta_{22}} &{\rm in}\,\,\R^N,
 \end{array}
 \right.
 \end{eqnarray*}
which contradicts Theorem 1.3 in \cite{FQ2} since (\ref{1.8}).

{\bf Case II: strongly coupled.} Choosing 
\[
\sigma_{1}=\frac{2\alpha(\beta_{12}+1)}{\beta_{12}\beta_{21}-1}\quad{\rm and}\quad \sigma_{2}=\frac{2\alpha(\beta_{21}+1)}{\beta_{12}\beta_{21}-1}.
\]
Since the system (\ref{fe}) is strongly coupled, we obtain
\begin{eqnarray}\nonumber
\sigma_1+2\alpha-\sigma_1\beta_{11}&>&0, \quad\sigma_1+2\alpha-\sigma_2\beta_{12}=0,\\\label{c}
\sigma_2+2\alpha-\sigma_1\beta_{21}&=&0, \quad\sigma_2+2\alpha-\sigma_2\beta_{22}>0.
\end{eqnarray}
and assumption $\gamma_{ij}<\beta_{ij}$  implies
\begin{eqnarray*}
\sigma_1+2\alpha-\sigma_1\gamma_{11}&>&0, \quad\sigma_1+2\alpha-\sigma_2\gamma_{12}>0,\\
\sigma_2+2\alpha-\sigma_1\gamma_{21}&>&0, \quad\sigma_2+2\alpha-\sigma_2\gamma_{22}>0.
\end{eqnarray*}
Using (A6), we have
\begin{eqnarray*}
\sigma_1+2\alpha-(\sigma_1+1)\theta_{11}&>&0, \quad\sigma_1+2\alpha-(\sigma_2+1)\theta_{12}>0,\\
\sigma_2+2\alpha-(\sigma_1+1)\theta_{21}&>&0, \quad\sigma_2+2\alpha-(\sigma_2+1)\theta_{22}>0.
\end{eqnarray*}
By a similar argument as before, we know the limit system of (\ref{fe}) is
\begin{eqnarray*}
\left\{\begin{array}{l@{\quad }l}
(-\Delta)^\alpha u\ge b(x_0)v^{\beta_{12}}\\
(-\Delta)^\alpha v\ge c(x_0)v^{\beta_{21}} 
 \end{array}
 \right.
 \end{eqnarray*}
in $\R^N$. Therefore, by Liouville type results  Theorem \ref{t5}, we come to a contradiction
as before if (\ref{1.9}) holds. We complete the prove. $\Box$\\

\begin{remark}\label{r}
We except Lemma \ref{l3.1} to hold in a large range such as (\ref{1.5}) and (\ref{1.6}). Unfortunately, this method of proof seems purely
 local and needs to be properly adapted to del with  non-local problems. Notice that there is no information for the functions $\tilde{u}_n$
 and $\tilde{v}_n$ in $\R^N\setminus \Om_n$, which leads difficult to pass limit appropriately in the system satisfies by $(\tilde{u}_n,\tilde{v}_n)$.
 \end{remark}

Now, we are in position to obtain a priori bounds of system (\ref{1.7}). We have already remarked that
due to the expected singularity of the gradient of solutions near boundary we need to work in spaces with weights which take care of the singularity. Hence, we fix $s$
satisfying
\beq\label{4.1}
0<s<1-\frac{\alpha}{\theta_{ij}}<1,
\eqq
where $\theta_{ij} (i,j=1,2)$ satisfying (A5) or (A6).
Let
\beq\label{4.3}
E=\{(u,v)\in C(\Om)\times C(\Om): \,\, \|u\|_1^{(-s)}<+\infty,\quad \|v\|_1^{(-s)}<+\infty\},
\eqq
where norm $\|\cdot\|_1^{(\tau)}$ is given  as in  (\ref{tau}) with $\tau=-s$.\\

Then we can prove

\begin{teo}\label{t4.1} Under the hypotheses of Theorem \ref{t2}, then there exists a constant $C>0$ such that for each couple of positive viscosity solutions $(u,v)$ of (\ref{1.7}) in 
$E$ with $s$ satisfying (\ref{4.1}), we have
\[
 \|u\|_1^{(-s)}, \|v\|_1^{(-s)}<C.
\]
\end{teo}

{\bf Proof.}  Assume that the conclusion of the theorem fails. Then there exists a sequence of positive solutions $(u_n,v_n)\in E$ such that 
$$\max\{ \|u_n\|_1^{(-s)},\|v_n\|_1^{(-s)}\}\rightarrow+\infty$$
as $n\rightarrow +\infty$. We may assume
\[
\lambda_n= \left(\|u_n\|_1^{(-s)}\right)^{-1/(\sigma_1+s)},
\]
if $\left(\|u_n\|_1^{(-s)}\right)^{\sigma_2+s}\ge \left(\|v_n\|_1^{(-s)}\right)^{\sigma_1+s}$ (up to a sequence), and
\[
\lambda_n=\left(\|v_n\|_1^{(-s)}\right)^{-1/(\sigma_2+s)}
\]
otherwise, for some constants $\sigma_1, \sigma_2>0$ which are to be determined later. Without loss of generality, we suppose that we are in the first of these two situations. 

Define
\[
M_n(x)=d(x)^{-s}u_n(x)+d(x)^{1-s}|\nabla u_n(x)|.
\]
Next, we choosing point $x_n\in \Om$ such that $M_n(x_n)\ge \sup_\Om M_n-\frac{1}{n}$ (the supremum may not be achieved). 
By our assumptions, we know $M_n(x_n)\rightarrow+\infty$. Let $\xi_n$ be a projection of $x_n$ on $\partial\Om$ and  let
\[
\tilde{u}_n(x)=\lambda_n^{\sigma_1}u_n(\xi_n+\lambda_n x)\quad{\rm and}\quad \tilde{v}_n(x)=\lambda_n^{\sigma_2} v_n(\xi_n+\lambda_n x)
\]
in $D_n=\{x\in \R^N:\,\,\xi_n+\lambda_nx\in \Om\}$.
Then function $(\tilde{u}_n, \tilde{v}_n)$ satisfies
\begin{eqnarray}\label{uee}
\left\{\begin{array}{l@{\quad }l}
(-\Delta)^\alpha \tilde{u}_n=\lambda_n^{\sigma_1+2\alpha-\sigma_1\beta_{11}}a(z)\tilde{u}_n^{\beta_{11}}+\lambda_n^{\sigma_1+2\alpha-\sigma_2\beta_{12}}b(z)\tilde{v}_n^{\beta_{12}}\\
\quad\quad+\lambda_n^{\sigma_1+2\alpha}h_1(z, \lambda_n^{-\sigma_1}\tilde{u}_n,\lambda_n^{-\sigma_2}\tilde{v}_n,\lambda_n^{-(\sigma_1+1)}\nabla\tilde{u}_n,\lambda_n^{-(\sigma_2+1)}\nabla\tilde{v}_n)&{\rm in}\,\,D_n,\\
(-\Delta)^\alpha \tilde{v}_n=\lambda_n^{\sigma_2+2\alpha-\sigma_1\beta_{21}}c(z)\tilde{u}_n^{\beta_{21}}+\lambda_n^{\sigma_2+2\alpha-\sigma_2\beta_{22}}b(z)\tilde{v}_n^{\beta_{22}}\\
\quad\quad+\lambda_n^{\sigma_2+2\alpha}h_2(z, \lambda_n^{-\sigma_1}\tilde{u}_n,\lambda_n^{-\sigma_2}\tilde{v}_n,\lambda_n^{-(\sigma_1+1)}\nabla\tilde{u}_n,\lambda_n^{-(\sigma_2+1)}\nabla\tilde{v}_n) &{\rm in}\,\,D_n,
 \end{array}
 \right.
 \end{eqnarray}
where $z=\lambda_n x+\xi_n$.
By assumption (A4) we have
\begin{eqnarray*}
|\lambda_n^{\sigma_1+2\alpha}h_1(z, \lambda_n^{-\sigma_1}\tilde{u}_n,\lambda_n^{-\sigma_2}\tilde{v}_n,\lambda_n^{-(\sigma_1+1)}\nabla\tilde{u}_n,\lambda_n^{-(\sigma_2+1)}\nabla\tilde{v}_n)|\\
\le C(\lambda_n^{\sigma_1+2\alpha}+\lambda_n^{\sigma_1+2\alpha-\sigma_1\gamma_{11}}|\tilde{u}_n|^{\gamma_{11}}+\lambda_n^{\sigma_1+2\alpha-\sigma_2\gamma_{12}}|\tilde{v}_n|^{\gamma_{12}}\\
+\lambda_n^{\sigma_1+2\alpha-(\sigma_1+1)\theta_{11}}|\nabla\tilde{u}_n|^{\theta_{11}}+\lambda_n^{\sigma_1+2\alpha-(\sigma_2+1)\theta_{12}}|\nabla\tilde{v}_n|^{\theta_{12}})
\end{eqnarray*}
and a similar estimate for $\lambda_n^{\sigma_2+2\alpha}h_2$.
Moreover, $(\tilde{u}_n,\tilde{v}_n)$ satisfying
\begin{eqnarray*}
\lambda_n^{s}d(\xi_n+\lambda_nx)^{-s}\tilde{u}_n(x)+\lambda_n^{s-1}d(\xi_n+\lambda_nx)^{1-s}|\nabla \tilde{u}_n(x)|
=\frac{M_n(\xi_n+\lambda_nx)}{M_n(x_n)}
\end{eqnarray*}
and 
\[
\lambda_n^{s}d(\xi_n+\lambda_nx)^{-s}\tilde{v}_n(x)\\+\lambda_n^{s-1}d(\xi_n+\lambda_nx)^{1-s}|\nabla \tilde{v}_n(x)|\le \lambda_n^{\sigma_2+t}\|v_n\|_1^{-s}.
\]
Then, using the fact $\lambda_n^{-1}d(\xi_n+\lambda_nx)=dist(x,\partial D_n)=:d_n(x)$ and the choice of the points $x_n$, we obtain for large $n$
\beq\label{4.4}
d_n(x)^{-s}\tilde{u}_n(x)+d_n(x)^{1-s}|\nabla\tilde{u}_n(x)|\le 2\quad{\rm in}\,\, D_n
\eqq
and 
\beq\label{4.44}
d_n(x)^{-s}\tilde{v}_n(x)+d_n(x)^{1-s}|\nabla\tilde{v}_n(x)|\le 1\quad{\rm in}\,\, D_n
\eqq
and moreover we know
\beq\label{4.5}
d_n(y_n)^{-s}\tilde{u}_n(y_n)+d_n(y_n)^{1-s}|\nabla\tilde{u}_n(y_n)|=1,
\eqq
where $y_n=\frac{x_n-\xi_n}{\lambda_n}$. 

Next, since $(u_n,v_n)$ solves system (\ref{1.1}), we can use Lemma \ref{l3.1} to obtain that 
\[
M_n(x_n)\le Cd(x_n)^{-s}(1+d(x_n)^{-\sigma_1})
\]
for some positive  constant $C$ independent of $n$. This implies that $d(x_n)\lambda_n^{-1}\le C$.
This bound entails that, passing to subsequence, $x_n\rightarrow x_0\in \partial \Om$ and $|y_k|=d(x_n)\lambda_n^{-1}\rightarrow d\ge0$ (in particular the points $\xi_n$ are uniquely determined at least for large $n$). Assuming that the outward unit normal to $\partial \Om$ at $x_0$ is $-e_N$, we also obtain that $D_n\rightarrow \R^N_+$ as $n\rightarrow +\infty$. 

We claim that $d>0$. To show this, notice that from (\ref{4.4}) and (\ref{4.44}) we have 
$$(-\Delta)^\alpha \tilde{u}_n\le C d_k^{\min\{(s-1)\theta_{11},(s-1)\theta_{12}\}}$$
 in $D_n$, for some constant independent of $n$. By our choice of $s$ and $\theta_{ij}\,\, (i,j=1,2)$, we have that
\beq\label{4.6}
s>\frac{\theta_{ij}-2\alpha}{\theta_{ij}}
\eqq
since $1<\theta_{ij}<2\alpha$.
Together with (\ref{4.1}), we see that
\beq\label{4.7}
\alpha<\max\{(1-s)\theta_{11},(1-s)\theta_{12}\}<2\alpha.
\eqq
Thus, by Lemma \ref{l2.3}, we can obtain that
\beq\label{4.8}
\tilde{u}_n(x)\le Cd_n(x)^{2\alpha+\min\{(s-1)\theta_{11},(s-1)\theta_{12}\}}, \quad{\rm when }\,\, d_n(x)<\delta.
\eqq
Furthermore, since $1<\theta_{ij}<2\alpha$ for $i,j=1,2$, 
\beq\label{4.10}
s>\frac{\theta_{ij}-2\alpha}{\theta_{ij}-1}.
\eqq 
Hence, $-s+2\alpha+(s-1)\theta_{ij}=s(\theta_{ij}-1)+2\alpha-\theta_{ij}>0$.
So, by (\ref{4.4}), we have
\[
\tilde{u}_n(x)\le 2d_n(x)^{s}\le \delta^{s-2\alpha-\min\{(s-1)\theta_{11},(s-1)\theta_{12}\}}d_n(x)^{2\alpha+\min\{(s-1)\theta_{11},(s-1)\theta_{12}\}}
\]
when $d_n(x)\ge \delta$. Hence $\|\tilde{u}_n\|_0^{(-2\alpha-\min\{(s-1)\theta_{11},(s-1)\theta_{12}\})} $ is bounded. Then we can use Lemma \ref{l2.2} with ($\tau=-\min\{(s-1)\theta_{11},(s-1)\theta_{12}\}$) to obtain
that
\beq\label{4.11}
|\nabla \tilde{u}_n(x)|\le Cd_n(x)^{2\alpha+\min\{(s-1)\theta_{11},(s-1)\theta_{12}\}-1}\quad{\rm in}\,\, D_n,
\eqq
where $C$ is also independent of $n$.  

Taking (\ref{4.8}) and (\ref{4.11}) in (\ref{4.5}), we deduce that
\[
1\le  Cd_n(y_n)^{2\alpha+\min\{(s-1)\theta_{11},(s-1)\theta_{12}\}-s}.
\]
 where $C$ is also independent of $n$.
This implies that $d_n(y_n)$ is bounded away from zero. Hence, $|y_n|$ is also since $0\in \partial D_n$. So that $d>0$ as claimed.

Finally, we can use Theorem \ref{2.14} together with  Ascoli-Arzel\'a's theorem and a diagonal argument to obtain that $\tilde{u}_n\rightarrow u$ and $\tilde{v}_n\rightarrow v$
in $C_{loc}^1(\R^N_+)$. Moreover, by (\ref{4.5}), we have 
$$d^{-s} u(y_0)+d^{1-s}|\nabla u(y_0)|=1$$
for some $y_0\in \R^N_+$. Hence, $(u,v)\not=(0,0)$ and $u(x)\le Cx_N^{2\alpha+\min\{(s-1)\theta_{11},(s-1)\theta_{12}} $ and $v(x)\le Cx_N^{2\alpha+\min\{(s-1)\theta_{12},(s-1)\theta_{22}} $ if $0<x_N<\delta$. 
Thus $u,v\in C(\R^N)$ and $u=v=0$ in $\R^N\setminus\R^N_+$. Passing to the limit in (\ref{uee}) with the help of Theorem \ref{2.13}, we obtain
\begin{eqnarray*}
\left\{\begin{array}{l@{\quad }l}
(-\Delta)^\alpha u=a(x_0)u^{\beta_{11}} &{\rm in}\,\,\R^N_+,\\
(-\Delta)^\alpha v=d(x_0)v^{\beta_{22}} &{\rm in}\,\,\R^N_+,\\
u=v=0 &{\rm in}\,\,\R^N\setminus \R^N_+
 \end{array}
 \right.
 \end{eqnarray*}
 in weakly coupled case and
 \begin{eqnarray*}
\left\{\begin{array}{l@{\quad }l}
(-\Delta)^\alpha u=b(x_0)v^{\beta_{12}} &{\rm in}\,\,\R^N_+,\\
(-\Delta)^\alpha v=c(x_0)u^{\beta_{21}} &{\rm in}\,\,\R^N_+,\\
u=v=0 &{\rm in}\,\,\R^N\setminus \R^N_+
 \end{array}
 \right.
 \end{eqnarray*}
 in strongly coupled case.
By Lemma \ref{l3.1}, we also know $u(x)\le C x_N^{-\sigma_1}$ and $v(x)\le Cx_N^{-\sigma_2}$ in $\R^N_+$ and thus $u$ and $v$ are bounded. By Theorems \ref{2.12} and \ref{2.9},
 we know $(u,v)$ is a classical solution. By the strong maximum principle, we have that $u>0$ and $v>0$ in $\R^N_+$. This is contradiction with Theorem 1.2 in \cite{QX} for scalar case and 
 Theorem \ref{t5} and Theorem 4.2 in \cite{QX} for system. We complete the proof.  $\Box$\\

\setcounter{equation}{0}
\section{ Existence results}

This section is devoted to prove  our existence results, Theorems  \ref{t1} and \ref{t2}. Both proofs are very similar, only that of Theorem \ref{t2} is slightly more involved. Hence, 
we only prove Theorem \ref{t2} here. 

 We assume $1/2<\alpha<1$. Fix $s$ verifying (\ref{4.1}) and consider the Banach space
$E$ (see (\ref{4.3})) with norm
\[
\|(u,v)\|_E=\max\{\|u\|_1^{(-s)},\|v\|_1^{(-s)}\},
\]
which is an ordered Banach space with the cone of nonnegative functions 
\[
K=\{(u,v)\in E:\,\, u\ge0\,\,{\rm and}\,\, v\ge0\,\, {\rm in}\,\,\Om\}.
\]
Observe that for every $(u,v)\in K$ we have
\begin{eqnarray*}
h_1(x,u,v,\nabla u,\nabla v)&\le& C d(x)^{\min\{(s-1)\beta_{11},(s-1)\beta_{12}\}}\\
h_2(x,u,v,\nabla u,\nabla v)&\le& C d(x)^{\min\{(s-1)\beta_{12},(s-1)\beta_{22}\}},
\end{eqnarray*}
where positive constant $C$ depending on the norms $\|u\|_1^{(-s)}$ and $\|v\|_1^{(-s)}$.
Moreover, as in the proof of Theorem \ref{t4.1}, we know
\[
\alpha<(s-1)\beta_{ij}<2\alpha, \quad i,j=1,2.
\]
Hence, applying Lemma \ref{l2.1} to system
\begin{eqnarray}\label{ss}
\left\{\begin{array}{l@{\quad }l}
(-\Delta)^\alpha u=a(x)\tilde{u}^{\beta_{11}}+b(x)\tilde{v}^{\beta_{12}}+h_1(x,\tilde{u},\tilde{v},\nabla\tilde{u},\nabla\tilde{v}) &{\rm in}\,\,\Omega,\\
(-\Delta)^\alpha v=c(x)\tilde{u}^{\beta_{12}}+d(x)\tilde{v}^{\beta_{22}}+h_2(x,\tilde{u},\tilde{v},\nabla\tilde{u},\nabla\tilde{v}) &{\rm in}\,\,\Omega,\\
u=v=0\,\,&{\rm in}\,\,\R^N\setminus\Omega,
 \end{array}
 \right.
 \end{eqnarray}
where $a,b,c,d$ satisfying (A1) and $h_1,h_2$ satisfying (A4),
then system (\ref{ss}) has a unique nonnegative solution $(u,v)$ with $\|u\|_0^{(-s)}<+\infty$ and $\|v\|_0^{(-s)}<+\infty$.Therefore, $(u,v)\in E$. So we define an operator $\Phi:\,\,K\rightarrow K$ by means of $(u,v)=\Phi(\tilde{u},\tilde{v})$. It is clear that nonnegative solutions of (\ref{1.1}) in $E$ coincide with the fixed points of the operator.\\

We first show the basic property of operator $\Phi$.

\begin{lemma}\label{l5.1}
The operator $\Phi:\,\,K\rightarrow K$ is compact.
\end{lemma}

{\bf Proof.}  We just need to do a slight modification of the proof of Lemma 11 in \cite{BDGQ} and thus we omit it here. $\Box$\\

The proof of our existence result is an application of degree theory for compact operators in cones. We start by recalling the  following well-known result (cf. Theorem 3.6.3 in \cite{K}).
\begin{teo}\label{t5.1}
Let $E$ is an ordered Banach space with positive cone $K$, and $U\subset K$ is an open bounded set containing $0$. Let $r>0$ such that
$B_r(0)\cap K\subset U$. Assume that $\Phi:\,\, U\rightarrow K$ is compact and satisfies\\
(i): for every $\mu\in [0,1)$, we have $u\not=\mu \Phi(u)$ for every $u\in K$ with $\|u\|=r$;\\
(ii): there exists $\phi\in K\setminus \{0\}$ such that $u-\Phi(u)=\rho\phi$, for every $u\in \partial U$ and every $\rho\ge0$. 

Then $\Phi$ has a fixed point in $U\setminus B_r(0)$.
\end{teo}

Now we are in a position  to prove Theorem \ref{t2}.\\

{\bf Proof of Theorem \ref{t2}.} We will show Theorem \ref{t5.1} is applicable to the operator $\Phi$ in $K\subset E$.

We first check first hypothesis $(i)$ in Theorem \ref{t5.1}. Assume we have $(u,v)=\mu \Phi(u,v)$ for some $\mu\in [0,1)$ and $(u,v)\in K$. This equivalent
\begin{eqnarray*}
\left\{\begin{array}{l@{\quad }l}
(-\Delta)^\alpha u=\mu(a(x)u^{\beta_{11}}+b(x)v^{\beta_{12}}+h_1(x,u,v,\nabla u,\nabla v) )&{\rm in}\,\,\Omega,\\
(-\Delta)^\alpha v=\mu(c(x)u^{\beta_{12}}+d(x)v^{\beta_{22}}+h_2(x,u,v,\nabla u,\nabla v))&{\rm in}\,\,\Omega,\\
u=v=0\,\,&{\rm in}\,\,\R^N\setminus\Omega.
 \end{array}
 \right.
 \end{eqnarray*}
By our  assumptions (A1) and (A4), we get that the right hand side of the above system can be bounded by 
\begin{center}
$\mu(a(x)u^{\beta_{11}}+b(x)v^{\beta_{12}}+h_1(x,u,v,\nabla u,\nabla v) )\le$\\
$Cd^{\min\{(s-1)\theta_{11},(s-1)\theta_{12}\}}\left(\|w\|_{E}^{\beta_{11}}+\|w\|_{E}^{\beta_{12}}+\|w\|_{E}^{\gamma_{
11}}+\|w\|_{E}^{\gamma_{12}}+\|w\|_{E}^{\theta_{11}}+\|w\|_{E}^{\theta_{12}}\right).$
\end{center}
where $w=(u,v)$ here and what follows and a similar estimate for $\mu(c(x)u^{\beta_{12}}+d(x)v^{\beta_{22}}+h_2(x,u,v,\nabla u,\nabla v))$.

Hence, by Lemmas \ref{l2.1} and \ref{l2.2} and $\alpha<(1-s)\theta_{ij}<2\alpha$ for $i,j=1,2$, we have 
\[
\|w\|_{E}\le C\Sigma_{i,j=1,2}\left(\|w\|_{E}^{\beta_{ij}}+\|w\|_{E}^{\beta_{ij}}+\|w\|_{E}^{\gamma_{
ij}}+\|w\|_{E}^{\gamma_{ij}}+\|w\|_{E}^{\theta_{ij}}+\|w\|_{E}^{\theta_{ij}}\right),
\]
here we have used the fact $s-2\alpha<(s-1)\theta_{ij}<0\,\,(i,j=1,2)$.
Since $\beta_{ij},\gamma_{ij},\theta_{ij}>1$ for $i,j=1,2$, this implies that  $\|w\|_{E}\ge r$ for some small positive $r>0$.
Thus, there are no positive solution of $(u,v)=\mu \Phi(u,v)$ if $\|(u,v)\|_E=r$ and $\mu\in (0,1)$, and $(i)$ follows. 

Next, we check $(ii)$. We take $\phi\in K$ to be the unique solution  (cf. Theorem 3.1 in \cite{FQ3}) of the problem
 \begin{eqnarray*}
\left\{\begin{array}{l@{\quad }l}
(-\Delta)^\alpha \phi=1&{\rm in}\,\,\Omega,\\
\phi=0\,\,&{\rm in}\,\,\R^N\setminus\Omega.
 \end{array}
 \right.
 \end{eqnarray*}
We claim that there are no solutions in $K$ of the equation $(u,v)-\Phi(u,v)=\rho(\phi,\phi)$ if $\rho$ is large enough. For this purpose, we note that
this equation equivalent to 
\begin{eqnarray}\label{5.3}
\left\{\begin{array}{l@{\quad }l}
(-\Delta)^\alpha u=a(x)u^{\beta_{11}}+b(x)v^{\beta_{12}}+h_1(x,u,v,\nabla u,\nabla v) +\rho&{\rm in}\,\,\Omega,\\
(-\Delta)^\alpha v=c(x)u^{\beta_{12}}+d(x)v^{\beta_{22}}+h_2(x,u,v,\nabla u,\nabla v)+\rho&{\rm in}\,\,\Omega,\\
u=v=0\,\,&{\rm in}\,\,\R^N\setminus\Omega.
 \end{array}
 \right.
 \end{eqnarray}
 Fix $\mu_1,\mu_2>\lambda_1$, where $\lambda_1$ is
\begin{center}
$ \lambda_1=\sup\{\lambda\in \R: {\rm there\,\, exists}\,\, u\in C(\R^N),\,\,u>0\,\, {\rm in}\,\, \Om\,\,
 {\rm with}\,\, u=0\,\,{\rm in} \,\,\R^N\setminus\Om\,\, {\rm such\,\,that}\,\,(-\Delta)^\alpha u\ge\lambda u\,\,{\rm in}\,\,\Om\}.$
\end{center}
By Lemma 13 in \cite{BDGQ}, we know $\lambda_1<+\infty$.

Since $\theta_{ij}>1$ for $i,j=1,2$ and $h_1,h_2$  are nonnegative, there exists positive constants $C_1$ and $C_2$ such that we have either (Case I)
\begin{eqnarray*}
(-\Delta)^\alpha u&\ge&\mu_1u-C_1+\rho,\\
(-\Delta)^\alpha v&\ge&\mu_2v-C_2+\rho,
\end{eqnarray*}
in $\Om$, or (Case II)
\begin{eqnarray*}
(-\Delta)^\alpha u&\ge&\mu_1v-C_1+\rho,\\
(-\Delta)^\alpha v&\ge&\mu_2u-C_2+\rho,
\end{eqnarray*}
in $\Om$. For Case I, we can take a similar argument as the proof in \cite{BDGQ} since $\mu_1,\mu_2>\lambda_1$ and choose $\rho\ge C_i\,\,(i=1,2)$. For Case II, we let $\rho\ge \max\{C_1,C_2\}/2$,  then 
\[
(-\Delta)^\alpha (u+v)\ge \min\{\mu_1,\mu_2\}(u+v) \quad{\rm in }\,\,\Om.\]
This contradicts the choice of $\mu_1,\mu_2$ and the definition of $\lambda_1$. Therefore, $\rho\le C$, and (\ref{5.3}) does not admit positive solutions in $E$ if $\rho$ is large.

Finally, since $h_1+\rho$ and $h_2+\rho$ also verifies (A4) for $\rho\le C$,  we can apply Theorem \ref{t5.1} to obtain that the solutions of (\ref{5.3}) are a priori bounded,
that is, there exists $R>r$ such that $\|(u,v)\|_E\le R$ for every positive solution of  (\ref{5.3}) with $\rho\ge 0$. Thus Theorem \ref{t5.1} is applicable with $U=B_R(0)\cap K$ and the
existence of a solution in $K$ follows. By the maximum principle, the solution is also positive. We complete the prove. $\Box$

\setcounter{equation}{0}
\section{ Acknowledgements}
A. Quaas was partially supported by Fondecyt Grant No. 1151180 Programa Basal, CMM. U. de Chile and Millennium Nucleus Center for Analysis of PDE NC130017.

\end{document}